# On Surfaces of finite Chen-type

Stylianos Stamatakis and Hassan Al-Zoubi
Aristotle University of Thessaloniki
Department of Mathematics
GR-54124 Thessaloniki, Greece
e-mail: stamata@math.auth.gr

**Abstract**

We investigate some relations concerning the first and the second Beltrami operators corresponding to the fundamental forms *I*, *II*, *III* of a surface in the Euclidean space $E^3$ and we study surfaces which are of finite type in the sense of B.-Y. Chen with respect to the fundamental forms *II* and *III*.

**2000 Mathematics Subject Classification**: 53A05, 53A10

**Keywords**: Surfaces in Euclidean space, Surfaces of finite type, Beltrami operator

**Introduction.** Let *M* be a (connected) surface in a Euclidean 3-space $E^3$ referred to any system of coordinates $u^1, u^2$, and let $n: M \to S^2 \subset E^3$ be its Gauss map. Denote by $\Delta^I$ the Laplace operator corresponding to the first fundamental form

$$I = g_{ij}\, du^i\, du^j.$$

Then the position vector $x = x(u^1, u^2)$ and the mean curvature *H* of *M* in $E^3$ satisfy the relation

(1) $\qquad \Delta^I x = -2H\, n,$

(with sign convention such that $\Delta^I = -\frac{\partial^2}{\partial x^2} - \frac{\partial^2}{\partial y^2}$ for the metric $ds^2 = dx^2 + dy^2$). From (1) we know the following two facts [7]

- $M$ is minimal if and only if all coordinate functions of $x$ are eigenfunctions of $\Delta^I$ with eigenvalue 0.

- $M$ lies in an ordinary sphere $S^2$ if and only if all coordinate functions of $x$ are eigenfunctions of $\Delta^I$ with a fixed nonzero eigenvalue.

In 1983 B.-Y. Chen introduced the notion of *Euclidean immersions of finite type* [1]. In terms of B.-Y. Chen theory a surface $M$ is said to be of finite type if its coordinate functions are a finite sum of eigenfunctions of its Laplacian $\Delta^I$. Then the above facts can be stated as follows

- $M$ is minimal if and only if $M$ is of null 1-type.

- $M$ lies in an ordinary sphere $S^2$ if and only if $M$ is of 1-type.

B.-Y. Chen posed in [2], [3] the problem of classifying the finite type surfaces in $E^3$. Many authors have been concerned with this problem. For a more detailed report and a recent survey we refer the reader to [3] and [4]. Up to this moment the only known finite Chen-type surfaces in $E^3$ are portions of spheres, circular cylinders and minimal surfaces.

In this paper we consider the second Beltrami operators $\Delta^{II}$ and $\Delta^{III}$ corresponding to the second and to the third fundamental form respectively of a surface $M$ in $E^3$. Our purpose is to find surfaces of finite Chen-type with respect to the second and to the third fundamental form. In paragraph 1 formulas for $\Delta^{II}x$, $\Delta^{III}x$, $\Delta^{II}n$ and $\Delta^{III}n$ are established. In paragraphs 2-4 we state and prove the main results.

Throughout this paper we make full use of the tensor calculus. The reader is referred to [5] for symbols and formulas.

**1.** We consider a surface $M$ which does not contain parabolic points and we denote by $b_{ij}$ the components of the second and by $e_{ij}$ the components of the third fundamental form of $M$. Furthermore, when we deal with the second fundamental form, we suppose that $M$ consists only of elliptic points. Let $\varphi$ and $\psi$ be two sufficient differentiable functions on $M$. Then the first differ-ential parameter of Beltrami with respect to the fundamental form $J = I, II, III$ of $M$ is defined by [5]

(1.1) $\quad\quad\quad \nabla^J(\varphi,\psi) := a^{ij}\, \varphi_{/i}\, \psi_{/j},$

where $\varphi_{/i} := \dfrac{\partial \varphi}{\partial u^i}$ and $(a^{ij})$ denotes the inverse tensor of $(g_{ij})$, $(b_{ij})$ and $(e_{ij})$ for $J = I, II$ and $III$ respectively.

We first prove the following relations:

(1.2) $\quad\quad\quad \nabla^I(\varphi,x) + \nabla^{II}(\varphi,n) = 0,$

(1.3) $\quad\quad\quad \nabla^{II}(\varphi,x) + \nabla^{III}(\varphi,n) = 0.$

For the proof of (1.2) we use (1.1) and the Weingarten equations [5, p.128]

(1.4) $\quad\quad\quad n_{/j} = -e_{jk}\, b^{km}\, x_{/m} = -b_{jk}\, g^{km}\, x_{/m},$

to obtain

$$\nabla^{II}(\varphi,n) = b^{ij}\, \varphi_{/i}\, n_{/j} = -b^{ij}\, \varphi_{/i}\, b_{jk}\, g^{km}\, x_{/m} = -g^{im}\, \varphi_{/i}\, x_{/m} = -\nabla^I(\varphi,x),$$

being (1.2).

We have similarly

$$\nabla^{III}(\varphi,n) = e^{ij}\, \varphi_{/i}\, n_{/j} = -e^{ij}\, \varphi_{/i}\, e_{jk}\, b^{km}\, x_{/m} = -b^{im}\, \varphi_{/i}\, x_{/m} = -\nabla^{II}(\varphi,x),$$

which is (1.3).

The second differential parameter of Beltrami with respect to the fundamental form $J = I, II, III$ of $M$ is defined by [5]

$$\Delta^J \varphi := -a^{ij}\, \nabla_i^J\, \varphi_{/j},$$

where $\varphi$ is a sufficient differentiable function, $\nabla_i^J$ is the covariant derivative in the $u^i$ direction with respect to the fundamental form $J$ and $(a^{ij})$ stands as in definition (1.1) for the inverse tensor of $(g_{ij})$, $(b_{ij})$ and $(e_{ij})$ for $J = I, II$ and $III$ respectively.

We first compute $\Delta^J x$ for $J = II$ and $III$.

Recalling the equations [5, p.128]

$$\nabla_j^{II} x_{/i} = -\frac{1}{2} b^{km} \nabla_k^I b_{ij} x_{/m} + b_{ij} n,$$

and inserting these into

$$\Delta^{II} x = -b^{ij} \nabla_j^{II} x_{/i}$$

we have

(1.5) $\qquad \Delta^{II} x = \frac{1}{2} b^{ij} b^{km} \nabla_k^I b_{ij} x_{/m} - b^{ij} b_{ij} n.$

By using the Mainardi-Codazzi equations [5, p.21]

(1.6) $\qquad \nabla_k^I b_{ij} - \nabla_i^I b_{jk} = 0,$

relation (1.5) becomes

(1.7) $\qquad \Delta^{II} x = \frac{1}{2} b^{ij} b^{mk} \nabla_i^I b_{jk} x_{/m} - 2 n.$

We consider the Christoffel symbols of the second kind with respect to the first and second fundamental form respectively

$$\Gamma_{ij}^{\ k} := \frac{1}{2} g^{km} (-g_{ij/m} + g_{im/j} + g_{jm/i}), \quad \Pi_{ij}^{\ k} := \frac{1}{2} b^{kr} (-b_{ij/m} + b_{im/j} + b_{jm/i}),$$

and we put

(1.8) $\qquad T_{ij}^{\ k} := \Gamma_{ij}^{\ k} - \Pi_{ij}^{\ k}.$

It is known that [5, p.125]

(1.9) $\qquad T_{ij}^{\ k} = -\frac{1}{2} b^{km} \nabla_m^I b_{ij},$

and

(1.10) $\qquad \Gamma_{ij}^{\ j} := \frac{g_{/i}}{2g}, \quad \Pi_{ij}^{\ j} := \frac{b_{/i}}{2b},$

where $g := \det(g_{ij})$ and $b := \det(b_{ij})$.

For the Gauss curvature $K$ of $M$ we know that

$$K = \frac{b}{g}.$$

From this, by using (1.8), (1.9) and (1.10), we find

(1.11) $$\frac{K_{/k}}{K} = \frac{b_{/k}}{b} - \frac{g_{/k}}{g} = 2(\Pi_{kj}{}^j - \Gamma_{kj}{}^j) = -2T_{kj}{}^j = b^{ij}\,\nabla^I_i\, b_{kj},$$

therefore

$$\frac{1}{2K}\,\nabla^{III}(K,n) = \frac{1}{2K}\,e^{ks}\, K_{/k}\, n_{/s} = \frac{1}{2}\,e^{ks}\, b^{ij}\,\nabla^I_i\, b_{kj}\, n_{/s}.$$

Hence, because of (1.4), we find

$$\frac{1}{2K}\,\nabla^{III}(K,n) = -\frac{1}{2}\,e^{ks}\, b^{ij}\, e_{sr}\, b^{rm}\,\nabla^I_i\, b_{kj}\, x_{/m} = -\frac{1}{2}\,b^{ij}\, b^{km}\,\nabla^I_i\, b_{kj}\, x_{/m}.$$

By simple substitution in (1.7) we obtain

(1.12) $$\Delta^{II}x = -\frac{1}{2K}\,\nabla^{III}(K,n) - 2\,n.$$

We use now the equations [5, p.128]

$$\nabla^{III}_j\, x_{/i} = -b^{km}\,\nabla^I_m\, b_{ij}\, x_{/k} + b_{ij}\, n$$

and

$$\Delta^{III}x = -e^{ij}\,\nabla^{III}_i\, x_{/j}$$

to get

(1.13) $$\Delta^{III}x = e^{ij}\, b^{km}\,\nabla^I_m\, b_{ij}\, x_{/k} - e^{ij}\, b_{ij}\, n.$$

We consider the Christoffel symbols of the second kind with respect to the third fundamental form

$$\Lambda_{ij}{}^k := \frac{1}{2}\,e^{km}\,(-e_{ij/m} + e_{im/j} + e_{jm/i}),$$

and we put

$$\widetilde{T}_{ij}{}^k := \Lambda_{ij}{}^k - \Pi_{ij}{}^k.$$

It is known that [5, p.22]

(1.14) $\quad \tilde{T}_{ij}{}^{k} = -\dfrac{1}{2} b^{km} \nabla^{III}_{m} b_{ij}$

and

(1.15) $\quad T_{ij}{}^{k} + \tilde{T}_{ij}{}^{k} = 0.$

On the other hand, using Ricci's Lemma

(1.16) $\quad \nabla^{III}_{j} e^{ik} = 0$

and the formula

(1.17) $\quad \dfrac{2H}{K} = e^{ik} b_{ik},$

we have

(1.18) $\quad \left(\dfrac{2H}{K}\right)_{/m} = \nabla^{III}_{m}(e^{ik} b_{ik}) = e^{ik} \nabla^{III}_{m} b_{ik}.$

By combining of (1.9), (1.14), (1.15) and (1.18) we obtain

$$e^{ij} b^{km} \nabla^{I}_{m} b_{ij} = -2 e^{ij} T_{ij}{}^{k} = 2 e^{ij} \tilde{T}_{ij}{}^{k} = -e^{ij} b^{km} \nabla^{III}_{m} b_{ij} = -b^{km} \left(\dfrac{2H}{K}\right)_{/m}$$

and so

(1.19) $\quad e^{ij} b^{km} \nabla^{I}_{m} b_{ij} x_{/k} = -b^{km} \left(\dfrac{2H}{K}\right)_{/m} x_{/k} = -\nabla^{II}\left(\dfrac{2H}{K}, x\right).$

From (1.13), (1.17) and (1.19) we find

$$\Delta^{III} x = -\nabla^{II}\left(\dfrac{2H}{K}, x\right) - \dfrac{2H}{K} n.$$

Finally, using (1.3) we arrive at

(1.20) $\quad \Delta^{III} x = \nabla^{III}\left(\dfrac{2H}{K}, n\right) - \dfrac{2H}{K} n.$

We focus now our interest to the computation of $\Delta^{J} n$ for $J = II$ and $III$. Firstly we mention the well-known formula

$$\Delta^{I} n = 2 \nabla^{I}(H, x) + 2(2H^2 - K) n.$$

Stamatakis and Al-Zoubi

Next we take into consideration the equations [5, p.128]

$$\nabla^{II}_i n_{/j} = -\frac{1}{2} b^{km} \nabla^{III}_m b_{ij} n_{/k} - e_{ij} n,$$

so that

$$\Delta^{II} n = -b^{ij} \nabla^{II}_i n_{/j}$$

takes the form

$$\Delta^{II} n = \frac{1}{2} b^{ij} b^{km} \nabla^{III}_m b_{ij} n_{/k} + b^{ij} e_{ij} n.$$

On account of

$$2H = b_{ik} g^{ik} = e_{ik} b^{ik},$$

and (1.14) we obtain

$$\Delta^{II} n = -b^{ij} \tilde{T}_{ij}^{\ k} n_{/k} + 2H n.$$

On use of (1.6), (1.9) and (1.15) we have

(1.21) $\qquad \Delta^{II} n = b^{km} T_{mj}^{\ j} n_{/k} + 2H n.$

On the other hand using (1.11) we have

$$\frac{1}{2K} \nabla^{II}(K,n) = \frac{1}{2K} b^{km} K_{/m} n_{/k} = -b^{km} T_{mj}^{\ j} n_{/k}.$$

Inserting this in (1.21) we get in view of (1.2)

(1.22) $\qquad \Delta^{II} n = \frac{1}{2K} \nabla^{I}(K,x) + 2H n.$

Finally from [5, p.128]

(1.23) $\qquad \nabla^{III}_k n_{/i} = -e_{ik} n$

we have

$$\Delta^{III} n = -e^{ik} \nabla^{III}_k n_{/i} = e^{ik} e_{ik} n,$$

so that we conclude

(1.24)    $\Delta^{III} n = 2\, n.$

**Remark**. The first-named author proved in [6] relations (1.20) and (1.24) using Cartan's method of the moving frame.

**2.** Outgoing from [1] we say that a surface $M$ is of *finite type* with respect to the fundamental form $J$, or briefly of *finite J-type*, $J = II, III$, if the position vector $x$ of $M$ can be written as a finite sum of nonconstant eigenvectors of the operator $\Delta^J$, that is if

(2.1)    $x = c + \sum_{i=1}^{m} x_i, \quad \Delta^J x_i = \lambda_i x_i, \quad i = 1,\ldots,m,$

where $c$ is a constant vector and $\lambda_1, \lambda_2, \ldots, \lambda_m$ are eigenvalues of $\Delta^J$; when there are exactly $k$ nonconstant eigenvectors $x_1, \ldots, x_k$ appearing in (2.1) which all belong to different eigenvalues $\lambda_1, \ldots, \lambda_k$, then $M$ is said to be of *J-type k*, and when $\lambda_i = 0$ for some $i = 1, \ldots, k$, then $M$ is said to be of *null J-type k*.

**Theorem 2.1.** *A surface $M$ in $E^3$ is of II-type* 1 *if and only if $M$ is part of a sphere.*

**Proof.** Let $M$ be a part of a sphere of radius $r$ centered at the origin. Then

(2.2)    $H = \dfrac{1}{r}, \quad K = \dfrac{1}{r^2}, \quad n = -\dfrac{1}{r} x.$

By using (1.12) and (2.2) we find

$$\Delta^{II} x = \dfrac{2}{r} x.$$

Therefore $M$ is of *II*-type 1 with corresponding eigenvalue $\lambda = \dfrac{2}{r}$.

Conversely, let $M$ be of *II*-type 1. Then we have

(2.3)    $\Delta^{II} x = \lambda x, \quad \lambda \in \mathbb{R}, \quad \lambda \neq 0.$

From (1.1), (1.12) and (2.3), we get

(2.4)    $2\lambda K\, x = -e^{ik} K_{/i}\, n_{/k} - 4K\, n,$

whence we find for the inner product of $x$ and $n$

(2.5) $$<x, n> = \frac{-2}{\lambda}.$$

Differentiating covariantly (2.4) in the $u^j$ direction with respect to the third fundamental form *III* and by using (1.23), (2.4) and Ricci's Lemma (1.16) we find

$$2\lambda K_{/j} x + 2\lambda K x_{/j} = -e^{ik} \nabla_j^{III} K_{/i} n_{/k} - 3K_{/j} n - 4K n_{/j}.$$

Taking the inner product of both sides of the last equation with *n*, we find in view of (2.5) $K_{/j} = 0$, $j = 1, 2$. Thus $K = const$. Then (2.4) yields $\lambda x = -2 n$, from which we have $|x| = \left|\frac{2}{\lambda}\right|$ and *M* is part of a sphere. □

**Theorem 2.2.** *The Gauss map of a surface M in $E^3$ is of II-type* 1 *if and only if M is part of a sphere.*

**Proof.** Let *M* be a part of a sphere of radius r centered at the origin. From (1.22) and (2.2) we get

$$\Delta^{II} n = \frac{2}{r} n,$$

and the Gauss map of *M* is of *II*-type1 with eigenvalue $\lambda = \frac{2}{r}$.

Conversely, let the Gauss map of *M* be of *II*-type1. Then we have

$$\Delta^{II} n = \lambda n, \quad \lambda \in \mathbb{R}, \quad \lambda \neq 0.$$

From this relation and (1.22) we obtain

$$\frac{1}{2K} \nabla^I(K,x) + (2H - \lambda) n = 0,$$

whence we deduce that $H = const. \neq 0$ and $K = const. \neq 0$, because $\nabla^I(K,x)$ is tangent to M, while *n* normal to *M*. So *M* is part of a sphere. □

**3.** We turn now to the study of surfaces, which are of finite type with respect to the third fundamental form considering first the minimal surfaces.

For the position vector of a minimal surface $M$ we find from (1.20)

(3.1)     $\Delta^{III} x = 0x,$

thus $M$ is of null $III$-type 1.

Conversely, let $M$ be of null $III$-type 1. From (1.20) and (3.1) it follows

$$\nabla^{III}(\frac{2H}{K}, n) - \frac{2H}{K} n = 0.$$

But $\nabla^{III}(\frac{2H}{K}, n)$ is tangent to $M$ while $n$ is normal to $M$. Therefore $\frac{2H}{K}$ vanishes and $M$ is minimal. So we have [6]

**Theorem 3.1.** *A surface $M$ in $E^3$ is of null III–type* 1 *if and only if $M$ is minimal.*

From (1.24) we get [6]

**Corollary 3.2.** *The Gauss map of every surface $M$ in $E^3$ is of III-type* 1. *The corresponding eigenvalue is $\lambda = 2$.*

Next we prove

**Theorem 3.3.** *A surface $M$ in $E^3$ is of III-type* 1 *if and only if $M$ is part of a sphere.*

**Proof.** Let $M$ be a part of a sphere of radius $r$ centered at the origin. By using (1.20) and (2.2) we find

$$\Delta^{III} x = 2x$$

which means that $M$ is of $III$-type 1.

Conversely, let $M$ be a surface of $III$-type 1. Then

(3.2)     $\Delta^{III} x = \lambda x, \quad \lambda \in \mathbb{R}, \quad \lambda \neq 0.$

From (1.1), (1.20) and (3.2) we get

(3.3) $\quad \lambda x = e^{ik} \left( \dfrac{2H}{K} \right)_{/i} n_{/k} - \dfrac{2H}{K} n.$

We differentiate now covariantly in the $u^j$ direction with respect to the third fundamental form *III*. Then we use (1.23) and Ricci's Lemma (1.16) to deduce

$$\lambda x_{/j} = e^{ik} \left[ \nabla_j^{III} \left( \dfrac{2H}{K} \right)_{/i} \right] n_{/k} - 2 \left( \dfrac{2H}{K} \right)_{/j} n - \dfrac{2H}{K} n_{/j},$$

whence, taking the inner product of both sides with $n$, we infer that $\left( \dfrac{2H}{K} \right)_{/j}$ = 0 for $j = 1,2$. Thus $\dfrac{H}{K} = \mu$ = const. Consequently (3.3) yields $\lambda x = -2\mu n$, which implies $|x| = \left| \dfrac{2\mu}{\lambda} \right|$ and $M$ is obviously part of a sphere. $\square$

**4.** Let $M$ be a surface of *III*-type $k$ with position vector $x = x(u^1, u^2)$ and $M^*$ be a parallel surface of $M$ in (directed) distance $\rho$ = *const.* $\neq 0$, so that

$$1 - 2\rho H + \rho^2 K \neq 0.$$

Then $M^*$ possesses the position vector

(4.1) $\quad x^* = x + \rho n.$

Denoting by $H^*$ the mean and by $K^*$ the Gauss curvature of $M^*$ we mention the following relation [5]

(4.2) $\quad \dfrac{H^*}{K^*} = \dfrac{H}{K} - \rho.$

Furthermore the surfaces $M$, $M^*$ have common unit normal vector and spherical image. Thus $III = III^*$ and

(4.3) $\quad \Delta^{III^*} = \Delta^{III}.$

For the surface $M$ we have

$$x = x_1 + x_2 + \ldots + x_k, \quad \Delta^{III} x_i = \lambda_i x_i, \quad i = 1, 2, \ldots, k.$$

On use of (1.24), (4.1) and (4.3) we observe that the position vector of the surface $M^*$ can be written as

$$x^* = x_1 + x_2 + \ldots + x_k + x_{k+1}, \quad x_{k+1} := \rho\, n,$$

where

$$\Delta^{III^*} x_i = \lambda_i x_i \quad \text{for} \quad i = 1, 2, \ldots, k,$$

and

$$\Delta^{III^*} x_{k+1} = 2\, x_{k+1}.$$

This result implies the following [6]

**Theorem 4.1.** *Every parallel surface $M^*$ of a surface $M$ of finite III-type $k$ is of $III^*$- type $k$ or $k+1$.*

This theorem includes obviously the following special case, which we express as

**Corollary 4.2.** *When $M$ is of null III-type $k$, then every parallel surface $M^*$ of $M$ is of null $III^*$- type $k$ or $k+1$.*

We prove now the following

**Theorem 4.3.** *Every parallel surface $M^*$ of a minimal $M$ is of null $III^*$-type 2.*

**Proof.** Let $M$: $x = x(u^1, u^2)$ be the minimal surface and $M^*$ the parallel in distance $\rho$, $\rho \in \mathbb{R} - \{0\}$. According to Corollary 4.2 $M^*$ is of null $III^*$-type 1 or 2. Let $M^*$ be of null $III^*$-type 1. Then by Theorem 3.1 $M^*$ is minimal. Thus there would be in the bundle of the parallel surfaces of $M$ two minimal surfaces, i.e. $M$ and $M^*$. This is possible only in case that $M$ is part of a plane. But then $K = 0$, which we have excluded.

So $M^*$ is of null $III^*$-type 2. □

Noticing from (4.2) that for the parallel $M^*$ of the minimal surface $M$ yields

Stamatakis and Al-Zoubi

$$\frac{H^*}{K^*} = -\rho = const. \neq 0,$$

we show the following

**Theorem 4.4.** *Let M be a surface satisfying* $\frac{H}{K} = \mu = const. \neq 0$, *which is not part of a sphere. Then M is of null III-type* 2.

**Proof.** We consider the parallel surface $M^*$ of $M$ in distance $\rho = \mu = \frac{H}{K}$. From (4.2) we have

$$\frac{H^*}{K^*} = 0,$$

consequently $M^*$ is minimal. But $M$ is the parallel of $M^*$ in distance $-\mu$. Therefore $M$, being parallel of a minimal surface, is of null *III*-type 2. □

**References**

[1] Chen, B.-Y.: Total mean Curvature and Submanifolds of Finite Type, World Scientific Publisher, 1984.

[2] Chen, B.-Y.: Surfaces of finite type in Euclidean 3-space, Bull. Soc. Math. Belg. Ser. B, **39** (1987), 243-254.

[3] Chen, B.-Y.: Some open problems and conjectures on submanifolds of finite type, Soochow J. Math., **17** (1991), 169-188.

[4] Chen, B.-Y.: A report on submanifolds of finite type, Soochow J. Math., **22** (1996), 117-337.

[5] Huck, H., Simon, U., Roitzsch, R., Vortisch, W., Walden, R., Wegner, B., Wendland, W.: Beweismethoden der Differentialgeometrie im Grossen, Lecture Notes in Mathematics. Vol. 335, 1973.

[6] Stamatakis, S.: Der 2. Beltramische Operator der dritten Grundform einer Fläche des E³, Proceedings of the 4[th] International Congress of Geometry, Thessaloniki (1996), 392-396.

[7] Takahashi, T.: Minimal immersions of Riemannian manifolds, J. Math. Soc. Japan **18** (1966), 380-385.